\documentclass[11pt,english]{amsart}
\usepackage[T1]{fontenc}
\usepackage[latin9]{inputenc}
\usepackage{amstext}
\usepackage{amsthm}
\usepackage{amssymb}

\makeatletter
\usepackage{amsmath,amsfonts,amssymb,amsthm,epsfig}

\usepackage{color}
\usepackage[final,allcolors=blue,colorlinks=true]{hyperref}

\def\R{\mathbb R}

\def\H{\mathbb H}                                      
\def\X{\mathbb X}

\def\al{\alpha}

\def\ep{\epsilon}

\def\ta{\theta}

\def\om{\omega}
\def\na{\nabla}
\def\Om{\Omega}  
\def\De{\Delta}      

\def\wq{\infty}
\def\pa{\partial}
\def\divergence{\text{\rm div}\,}

\newcommand{\medint}{-\kern -,375cm\int}         
\newcommand{\medintinrigo}{-\kern -,315cm\int}

 \newcommand{\dive}{{\rm div}}

\numberwithin{equation}{section}
\newtheorem{theorem}{Theorem}[section]
\newtheorem*{theorem*}{Theorem}  

\newtheorem*{expectation*}{Rivi\`ere's expectation}

\newtheorem*{conclusion*}{Conclusin}

\newtheorem*{conjecture*}{Conjecture}
\newtheorem{corollary}[theorem]{Corollary}
\newtheorem*{corollary*}{Corollary}

\newtheorem*{lemma*}{Lemma}

\newtheorem*{notation*}{Notation}

\newtheorem*{problem*}{Problem}

\newtheorem*{proposition*}{Proposition}

\newtheorem*{remark*}{Remark}

\newtheorem*{example*}{Example}
\newtheorem*{Acknowledgements*}{Acknowledgements}           
\theoremstyle{definition}

\makeatother

\usepackage{babel}
\begin{document}
\title[]{Global conservation law for an even order elliptic system with antisymmetric potential} 

    \author[C.-Y. Guo, C.-L. Xiang and  G.-F. Zheng ]{Chang-Yu Guo, Chang-Lin Xiang* and Gao-Feng Zheng}         

\address[Chang-Yu Guo]{Research Center for Mathematics and Interdisciplinary Sciences, Shandong University 266237,  Qingdao, P. R. China} 
\email{changyu.guo@sdu.edu.cn}

\address[Chang-Lin Xiang]{Three Gorges Mathematical Research Center, China Three Gorges University,  443002, Yichang,  P. R. China}
\email{changlin.xiang@ctgu.edu.cn}

\address[Gao-Feng Zheng]{School of Mathematics and Statistics, Central China Normal University, Wuhan 430079,  P. R.  China} 
\email{gfzheng@mail.ccnu.edu.cn}

\thanks{*Corresponding author}

\thanks{C.-Y. Guo is supported by the Young Scientist Program of the Ministry of Science and Technology of China (No.~2021YFA1002200), the National Natural Science Foundation of China (No.~12101362) and the Natural Science Foundation of Shandong Province (No.~ZR2021QA003). The corresponding author C.-L. Xiang is financially supported by the National Natural Science Foundation of China (No.~11701045). G.-F. Zheng is supported by the National Natural Science Foundation of China (No. 11571131). }

\begin{abstract}
In this note, we refine the local conservation law obtained by Lamm-Rivi\`ere for fourth order systems and de Longueville-Gastel for general even order systems to a global conservation law. 

\end{abstract}

\maketitle
{\small    
	\smallskip
	\newline
}
\bigskip

\section{Introduction}

Interesting geometric partial differential equations are usually of critical or even supercritical  nonlinearity in nature.  For instance, consider the harmonic mappings equation
\begin{equation}\label{eq:harmonic map}
	-\Delta u=A(u)(\nabla u,\nabla u), \qquad u\in W^{1,2}(B^m, N),
\end{equation}
where $N$ is a closed Riemannian manifold embedded isometrically in some $\R^n$ and $A$ dentoes the second fundamental form  $N$ in $ \R^n$.  Two fundamental analysis issues regarding \eqref{eq:harmonic map} are the regularity of  weak solutions and the weak sequential compactness of solutions with uniformly bounded energy. However, the $L^1$ integrability of the right hand side of \eqref{eq:harmonic map} makes it on the borderline of the usual elliptic regularity theory. In supercritical dimensions, that is when $m\geq 3$, weakly harmonic mappings  even could be everywhere discontinuous; see the surprising examples constructed in \cite{Riviere-1995-Acta}. Another equation with similar regularity difficulties is the so called prescribed mean curvature equation 
\begin{equation}\label{eq:prescribed MCE}
-\Delta u = -2H(u)u_x\wedge u_y, \qquad u\in W^{1,2}(B^2, \R^3),
\end{equation} where $H\in L^{\wq}(\R^3)$.


On the other hand, the lack of maximum principle for  elliptic systems makes the regularity issues  rather difficult problems. To solve both analysis issues mentioned above, the approach of conservation law turns out to be a very precious  tool that is remained to be effictive. Suppose for simplicity $N=\mathbb{S}^n$ is the sphere. Shatah \cite{Shatah-1988-CPAM} proved that $u$ is a solution of \eqref{eq:harmonic map} if and only if it satisfies the following conservation law
\begin{equation}\label{eq:conservation law for sphere}
	\dive(u^i\nabla u^j-u^j\nabla u^i)=0\qquad \text{for all }i,j\in \{1,\cdots,n\}.
\end{equation}
Combining \eqref{eq:conservation law for sphere} with the Rellich-Kondrachov compactness theorem leads rather directly to the affirmative answer of the weak compactness question.  In the conformal dimension $m=2$,   the answer to the regularity issue was settled by H\'elein using again \eqref{eq:conservation law for sphere}: by  Poincar\'e's lemma, there exists $B\in W^{1,2}$ such that $\nabla^\perp B_{ij}=u_i\nabla u_j-u_j\nabla u_i$, where $\nabla^\perp=(-\partial_y,\partial_x)$. Thus \eqref{eq:harmonic map} is equivalent to 
\begin{equation}\label{eq:Helein method}
	-\Delta u=\nabla^\perp B\cdot \nabla u.
\end{equation} 
Now the  curl-grad structure in the right-hand side of \eqref{eq:Helein method} implies better  regularity than being simply in $L^1$, which was a key observation discovered by Wente \cite{Went-1969} when studying constant mean curvature equations,  see also \cite{CLMS-1993} using the language of Hardy spaces. Using Wente's Lemma, one easily concludes that $u\in C^0(B^m,\mathbb{S}^n)$. For general target spaces other than $\mathbb{S}^n$, H\'elein \cite{Helein-2002} also solved in dimension two  the regularity issue via his famous moving frame method; based on the idea of conservation law but permits to avoid a direct conservation law. However, it seems that the approach of H\`elein does not apply to general conformally invariant  second order elliptic variational problems with quadratic growth; in particular, does not apply to the prescribed mean curvature equation \eqref{eq:prescribed MCE} under the weakest assumption that  $H\in L^{\wq}(\R^3)$. 

A  direct conservation law was not found until the  significant work \cite{Riviere-2007} of  Rivi\`ere, where
the author found  that not only the  harmonic mappings equation \eqref{eq:harmonic map}, but also the  prescribed mean curvature equation \eqref{eq:prescribed MCE} and the  Euler-Lagrange equation of the general  conformally invariant  second order elliptic  variational problems with quadratic growth in dimension two. More precisely,  he introduced  the  second order linear elliptic  system 
\begin{equation}\label{eq:Riviere 2007}
	-\Delta u=\Omega\cdot \nabla u \qquad \text{in }B^2,
\end{equation}
where $u\in W^{1,2}(B^2, \R^n)$ and $\Omega=(\Omega_{ij})\in L^2(B^2,so_n\otimes \Lambda^1\R^2)$. One can verify (see \cite{Riviere-2007}) that \eqref{eq:Riviere 2007} includes
the Euler-Lagrange equations of  critical points of all second order conformally invariant variational functionals which  act on mappings $u\in W^{1,2}(B^2,N)$ from  $B^2\subset \R^2$ into a closed Riemannian manifold $N\subset \R^n$. In particular, \eqref{eq:Riviere 2007} includes the equations of weakly harmonic mapping equation \eqref{eq:harmonic map} and the  prescribed mean curvature equation \eqref{eq:prescribed MCE}.  Since $\Omega\in L^2$, system \eqref{eq:Riviere 2007} is critical in the sense that $\Om\cdot \na u\in L^1(B^2)$,  which allows for discontinuous weak solutions in general. Due to antisymmetry of $\Om $, the application of (an adapted version of) Uhlenbeck's gauge theory \cite{Uhlenbeck-1982} allows  Rivi\`ere to write 
\begin{equation*}\label{eq:decompse Omega}
\Om=P^{-1}dP+P^{-1}d^{\ast}\xi P,
\end{equation*}
where $P\in W^{1,2}(B^2,SO(n))$ and $\xi\in W^{1,2}(B^2,so_n\otimes \wedge^2\R^2)$.
Then  Rivi\`ere succeeded in finding  functions $A\in L^{\wq}\cap W^{1,2}(B^2, Gl(n))$ and $B\in W^{1,2}(B^2, M_n)$ such that
$$\nabla A-A\Omega=\na^{\bot}B.$$
Once such $A$ and $B$ were found, it is straightforward to check that system \eqref{eq:Riviere 2007} can be written  equivalently as the direct conservation law
\begin{equation}\label{eq:conservation law of Riviere}
	\divergence(A\nabla u+B\nabla^{\bot} u )=0,
\end{equation}
from which everywhere continuity of  weak solutions of system \eqref{eq:Riviere 2007} can be derived. As applications, this recovered the famous regularity result of  H\'elein  \cite{Helein-2002}, and confirmed affirmatively two long-standing regularity conjectures by Hildebrandt and Heinz on  conformally invariant geometrical problems and the prescribed bounded mean curvature equations respectively, see \cite{Riviere-2007} for details.

Motivated by problems from conformal geometry, to search  higher order conformally invariant  mappings in higher dimensions that play a similar role as that of  harmonic mapping in dimension two, it is natural to consider the $m$-polyharmonic energy functionals
\begin{equation*}\label{eq:def polyharmonic energy}
	E_m(u):=\frac{1}{2}\int_{B^{2m}}|D^mu|^2dx, \qquad u\in W^{m,2}(B^{2m},N),
\end{equation*} 
or 
\begin{equation*}\label{eq: intrinsic polyharmonic energy}
I_m(u):=\frac{1}{2}\int_{B^{2m}}|\na^{m-1}Du|^2dx, \qquad u\in W^{m,2}(B^{2m},N),
\end{equation*} 
where $\nabla$ is the Levi-Civita connection on $N$ and  $m\ge 2$.  Critical points of $E_m$ ($I_m$ resp.) are called extrinsically (intrinsically resp.) $m$-polyharmonic mappings. When $m=2$, critical points of $E_2$ ($I_2$ resp.) are called extrinsically (intrinsically resp.) biharmonic mappings. 

Historically, Chang, Wang and Yang \cite{Chang-W-Y-1999} established a regularity theory for biharmonic mappings into the sphere and later Wang~\cite{Wang-2004-MZ,Wang-2004-CPAM} extended this regularity theory to biharmonic mappings into closed Riemannian manifolds. Inspired by Rivi\`ere's direct conservation law approach~\cite{Riviere-2007}  for the second order elliptic system \eqref{eq:Riviere 2007},  Lamm and Rivi\`ere proposed in \cite{Lamm-Riviere-2008} the fourth order elliptic system
\begin{equation}\label{eq:Lamm-Riviere 2008}
	\De^{2}u=\De(V\cdot\na u)+{\rm div}(w\na u)+W\cdot\na u \quad  \text{in }B^4,
\end{equation}
where $$V\in W^{1,2}(B^4,M_n\otimes \Lambda^1\R^{4}), w\in L^{2}(B^4,M_n)$$
and $W$ is of the form
$W=\na\om+F$
with $$\om\in L^{2}(B^4,so_n)$$ and $F\in L^{\frac{4}{3},1}(B^4,M_n\otimes \Lambda^1\R^{4})$. They verified that system \eqref{eq:Lamm-Riviere 2008} indeed includes both extrinsic and intrinsic biharmonic mappings from $B^4$ into closed Riemannian manifolds as special cases.
Similar to the second order case \cite{Riviere-2007}, if we were able  to find $A\in W^{2,2}\cap L^{\wq}(B^4,M_n)$ and $B\in W^{1,4/3}(B^4,M_n\otimes\wedge^{2}\R^{4})$
satisfying   
\begin{equation}\label{eq:fourth order for CL}
	d\Delta A+\Delta AV-\nabla Aw+AW=d^*B\qquad \text{ in }B^4,
\end{equation}
 then it is immediate  to check that $u$ is a solution of \eqref{eq:Lamm-Riviere 2008} if and only if it satisfies the direct conservation law 
\begin{equation}\label{eq:conservation law of Lamm Riviere}
	\operatorname{div}[\nabla(A \Delta u)-2 \nabla A \Delta u+\Delta A \nabla u-A w \nabla u+\nabla A(V \cdot \nabla u)-A \nabla(V \cdot \nabla u)-B \cdot \nabla u]=0
\end{equation} in $B^4$. 
Due to some technical difficulties, Lamm and Rivi\`ere~\cite{Lamm-Riviere-2008} only succeeded in finding $A\in W^{2,2}\cap L^{\wq}(B_{1/2}^4,M_n)$ and $B\in W^{1,4/3}(B_{1/2}^4,M_n\otimes\wedge^{2}\R^{4})$ such that \eqref{eq:fourth order for CL} holds in the smaller ball $B_{1/2}^4\subset B^4$. Consequently, their conservation law is ``local'' in the sense that it only holds in a strictly smaller region $B_{1/2}^4$ rather than on the whole domain $B^4$, where system \eqref{eq:Lamm-Riviere 2008} is defined.

The study on general $m$-polyharmonic mappings have also attracted great attention in the last decades. For instance, Gastel and Scheven \cite{Gastel-Scheven-2009CAG} obtained a regularity theory for both extrinsic and intrinsic polyharmonic mappings in critical dimensions via the moving frame method of Helein \cite{Helein-2002}. For more progress in this respect, see e.g. \cite{Goldstein-Strzelecki-Zatorska-2009,Lamm-Wang-2009} and the references therein. Based on the works of Rivi\`ere\cite{Riviere-2007} for $m=1$ and Lamm-Rivi\`ere \cite{Lamm-Riviere-2008} for $m=2$, it is natural to find a unified treatment for general $m$-polyharmonic mappings via the direct conservation law approach. 

This was achieved very recently by de Longueville and Gastel in their interesting work \cite{deLongueville-Gastel-2019}. To includes (both extrinsic and intrinsic) $m$-polyharmonic mappings, they introduced the  even order linear elliptic system
\begin{equation}\label{eq:Longue-Gastel system}
	\Delta^{m}u=\sum_{l=0}^{m-1}\Delta^{l}\left\langle V_{l},du\right\rangle +\sum_{l=0}^{m-2}\Delta^{l}\delta\left(w_{l}du\right) \qquad \text{ in }B^{2m}.
\end{equation}
System \eqref{eq:Longue-Gastel system} reduces to \eqref{eq:Lamm-Riviere 2008} when $m=2$, and to \eqref{eq:Riviere 2007} when $m=1$. The coefficient functions are assumed to satisfy
\begin{equation}\label{eq:coefficient w V}
	\begin{aligned}
		&w_{k} \in W^{2 k+2-m, 2}\left(B^{2 m}, \mathbb{R}^{n \times n}\right) \quad \text { for } k \in\{0, \ldots, m-2\} \\
		&V_{k} \in W^{2 k+1-m, 2}\left(B^{2 m}, \mathbb{R}^{n \times n} \otimes \wedge^{1} \mathbb{R}^{2 m}\right) \quad \text { for } k \in\{0, \ldots, m-1\}.
	\end{aligned}
\end{equation}
Moreover,  the first order potential $V_0$ has the decomposition $V_{0}=d \eta+F$ with
\begin{equation}\label{eq:coefficient eta F}
	\eta \in W^{2-m, 2}\left(B^{2 m}, s o(n)\right), \quad F \in W^{2-m, \frac{2 m}{m+1}, 1}\left(B^{2 m}, \mathbb{R}^{n \times n} \otimes \wedge^{1} \mathbb{R}^{2 m}\right).
\end{equation}

To formulate the conservation law of  de Longueville and Gastel \cite{deLongueville-Gastel-2019}, we set
\begin{equation}\label{eq:theta for small coefficient}
	\begin{aligned}
		\theta_{D}:=\sum_{k=0}^{m-2}&\|w_k\|_{W^{2k+2-m,2}(D)}+\sum_{k=1}^{m-1}\|V_k\|_{W^{2k+1-m,2}(D)}\\
		&+\|\eta\|_{W^{2-m,2}(D)}+\|F\|_{W^{2-m,\frac{2m}{m+1},1}(D)}
	\end{aligned}
\end{equation}
for $D\subset \R^{2m}$. Then, under the smallness assumption
\begin{equation}\label{eq:smallness assumption}
	\theta_{B^{2m}}<\ep_m,
\end{equation}
they were able to find $A\in W^{m,2}\cap L^\infty(B_{1/2}^{2m},Gl(n))$ and $B\in W^{2-m,2}(B^{2m}_{1/2},\R^{n\times n}\otimes \wedge^2\R^{2m})$  such that
\begin{equation}\label{eq:for CL general order}
	\Delta^{m-1}dA+\sum_{k=0}^{m-1}(\Delta^k A)V_k-\sum_{k=0}^{m-2}(\Delta^k dA)w_k=\delta B\qquad \text{in }B_{1/2}^{2m}.
\end{equation}
Consequently, they obtained the local conservation law:  $u$ solves \eqref{eq:Longue-Gastel system} in $B_{1/2}^{2m}$, if and only if it is a distributional solution of 
\begin{equation}\label{eq:conservation law of D-G}
	\begin{aligned}
		0&=\delta\Big[\sum_{l=0}^{m-1}\left(\Delta^{l} A\right) \Delta^{m-l-1} d u-\sum_{l=0}^{m-2}\left(d \Delta^{l} A\right) \Delta^{m-l-1} u \\ &\qquad -\sum_{k=0}^{m-1} \sum_{l=0}^{k-1}\left(\Delta^{l} A\right) \Delta^{k-l-1} d\left\langle V_{k}, d u\right\rangle+\sum_{k=0}^{m-1} \sum_{l=0}^{k-1}\left(d \Delta^{l} A\right) \Delta^{k-l-1}\left\langle V_{k}, d u\right\rangle \\ &\qquad -\sum_{k=0}^{m-2} \sum_{l=0}^{k-2}\left(\Delta^{l} A\right) d \Delta^{k-l-1} \delta\left(w_{k} d u\right)+\sum_{k=0}^{m-2} \sum_{l=0}^{k-2}\left(d \Delta^{l} A\right) \Delta^{k-l-1} \delta\left(w_{k} d u\right) \\ &\qquad -\langle B, d u\rangle\Big],
	\end{aligned}
\end{equation}
in $B_{1/2}^{2m}$, where $d \Delta^{-1} \delta$  denotes the identity map. We  mention that a similar local conservation law but with a slightly different form of $A$ is deduced by H\"orter and Lamm \cite{Horter-Lamm-2020}.

Comparing with the "global" conservation law \eqref{eq:conservation law of Riviere} of Rivi\`ere for the second oredr probelm \eqref{eq:Riviere 2007}, it is natural to formulate an open problem as follows.
\medskip 

\textbf{Problem A.} \emph{
	Can we establish a global conservation law for \eqref{eq:Longue-Gastel system}? That is, can we find $A\in W^{m,2}\cap L^\infty(B^{2m},Gl(n))$ and $B\in W^{2-m,2}(B^{2m},\R^{n\times n}\otimes \wedge^2\R^{2m})$ in the whole domain $B^{2m}$ such that \eqref{eq:for CL general order} and thus \eqref{eq:conservation law of D-G}  holds in $B^{2m}$.}
\medskip 

The aim of this note is to  give an affirmative answer to the above Problem. Our main result reads as follows.
\begin{theorem}[Global conservation law]\label{thm:general even order} 
	 There exist  constants $\ep_{m},C_m>0$ such that under the smallness assumption \eqref{eq:smallness assumption}, there exist $A\in W^{m,2}\cap L^\infty(B^{2m},Gl(n))$ and $B\in W^{2-m,2}(B^{2m},\R^{n\times n}\otimes \wedge^2\R^{2m})$ satisfying \eqref{eq:for CL general order}  in $B^{2m}$. Moreover, 
	 	\[
	 \|A\|_{W^{m,2}(B^{2m})}+\|dist(A,SO(m))\|_{L^{\infty}(B^{2m})}+\|B\|_{W^{2-m,2}(B^{2m})}\leq C_{m}\ta_{B^{2m}}.
	 \]
	 Consequently, $u$ solves \eqref{eq:Longue-Gastel system} if and only if it satisfies the conservation law \eqref{eq:conservation law of D-G}  on $B^{2m}$.  
\end{theorem}

We remark that it is able to find, by comining the proof of de Longueville and Gastel \cite{deLongueville-Gastel-2019} with the proof of Theorem \ref{thm:general even order} below, that $A\in W^{m+1, \frac{2m}{m+1}, 1}(B^{2m}, Gl(m)) $ and $B\in W^{3-m, \frac{2m}{m+1}, 1}(B^{2m}, M(m)) $ with the same estimate, which maybe have potential applications in future research.  We  remark that it is also  possible to deduce a global conservation law in the form of H\"orter and Lamm \cite{Horter-Lamm-2020}. We leave it to interested readers.

We separate the case $m=2$ from Theorem \ref{thm:general even order} as a corollary, since it improves \cite[Theorem 1.5]{Lamm-Riviere-2008}.
\begin{corollary}\label{coro:fourth case} 
There exist $\ep_{m},C_{m}>0$ satisfying the following property. For any $V,w,F,\om$ satisfying the smallness condition
	\[
	\|V\|_{W^{1,2}}+\|w\|_{L^{2}}+\|\om\|_{L^{2}}+\|F\|_{L^{4/3,1}}<\ep_{m},
	\]
	there exist $A\in W^{2,2}\cap L^{\wq}(B^{4},Gl(m))$ and $B\in W^{1,4/3}(B^{4},M_{m}\otimes\wedge^{2}\R^{4})$
	such that \eqref{eq:fourth order for CL} holds on $B^{4}$ with the following estimate
	and 
	\[
	\|A\|_{W^{2,2}}+\|dist(A,SO(m))\|_{L^{\infty}}+\|B\|_{W^{1,4/3}}\leq C_{m}\left(\|V\|_{W^{1,2}}+\|w\|_{L^{2}}+\|\Om\|_{L^{2}}+\|F\|_{L^{4/3,1}}\right).
	\]
\end{corollary}

It will be clear from the arguments in the next section that, there exist infinitely many $A,B$,  such that the conservation law of Theorem  \ref{thm:general even order} and Corollary \ref{coro:fourth case}  hold. Thus an interesting problem is to  classify all  $A,B$ so as to make it clear the essence of the conservation law. 

As the first and also the most fundamental application of conservation law, we would like to point out that once the conservation law \eqref{eq:conservation law of D-G} is derivated, everywhere continuity of weak solutions to \eqref{eq:Longue-Gastel system} follows easily from the Riesz potential theory and the Lorentz-Sobolev embedding  
\[W^{m,2,1}(B^{2m})\subset C^0(\overline{B^{2m}}).\] However, to improve the everywhere continuity to H\"older continuity, extra efforts are necessary. For the case $m=1$, H\"older continuity was obtained in \cite{Riviere-2012}; for the case $m=2$, see \cite{Guo-Xiang-2019-Boundary} and for general $m\in \mathbb{N}$, see \cite{Guo-Xiang-Zheng-2021-Lp}. Furthermore, in \cite{Sharp-Topping-2013-TAMS,Guo-Xiang-Zheng-2020-Lp,Guo-Xiang-Zheng-2021-Lp}, more than the H\"older continuity, an optimal $L^p$ theory for weak solutions to (inhomogenous version of) \eqref{eq:Riviere 2007}, \eqref{eq:Lamm-Riviere 2008} and \eqref{eq:Longue-Gastel system} was established. There are also alternative ways to derive H\"older continuity (or $L^p$ theory) of \eqref{eq:Riviere 2007}, \eqref{eq:Lamm-Riviere 2008} and \eqref{eq:Longue-Gastel system}; we refer the readers to \cite{Riviere-Struve-2008,Struwe-2008,Guo-Xiang-2019-Higher,Guo-Wang-Xiang-2021-arXiv} for more information.

The general idea for the proof of Theorem \ref{thm:general even order} is due to Rivi\`ere \cite{Riviere-2007} and thus is rather similar to the one presented in \cite{Lamm-Riviere-2008} and \cite{deLongueville-Gastel-2019}. We shall give a detailed treatment for the fourth order system \eqref{eq:Lamm-Riviere 2008} and then indicate the necessary changes that lead to the proof of Theorem \ref{thm:general even order}. The main new ingredient of our approach is to use an extension argument, instead of a cutoff argument used in \cite{Lamm-Riviere-2008,deLongueville-Gastel-2019}, to solve \eqref{eq:fourth order for CL} (or \eqref{eq:for CL general order})  in the enlarged region $B_2^{4}$ (or $B_2^{2m}$). A similar fix point argument as that of \cite{Lamm-Riviere-2008,deLongueville-Gastel-2019} gives us a solution $(A,B)$ of \eqref{eq:fourth order for CL} (or \eqref{eq:for CL general order})  in $B_2^{4}$ (or $B_2^{2m}$), which is also a solution to \eqref{eq:fourth order for CL} (or \eqref{eq:for CL general order})  in $B^{4}$ (or $B^{2m}$). 

As another fundamental  application of the conservation law in Theorem \ref{thm:general even order}, we deduce the following Palais-Smale sequential compactness for \eqref{eq:Longue-Gastel system}, which is expected to hold in view of the  second and fourth order case, see Rivi\`ere \cite{Riviere-2007} and Guo and Xiang \cite{Guo-Xiang-2019-Boundary}.
\begin{theorem}\label{thm:weak compactness}
	Let $f_i$ be a sequence in $H^{-m}:=(W^{m,2}(B^{2m},\R^n))^\ast$ that converges to 0 in $H^{-m}$ and $u_i$ be a bounded sequence in $W^{m,2}(B^{2m},\R^n)$ solving
	\begin{equation}\label{eq:seq compact}
	\Delta^{m}u_i=\sum_{l=0}^{m-1}\Delta^{l}\left\langle V_{l,i},du_i\right\rangle +\sum_{l=0}^{m-2}\Delta^{l}\delta\left(w_{l,i}du_i\right)+f_i \qquad \text{ in }B^{2m},
	\end{equation}
where
\begin{equation*}
	\begin{aligned}
		&w_{k,i} \in W^{2 k+2-m, 2}\left(B^{2 m}, \mathbb{R}^{n \times n}\right) \quad \text { for } k \in\{0, \ldots, m-2\} \\
		&V_{k,i} \in W^{2 k+1-m, 2}\left(B^{2 m}, \mathbb{R}^{n \times n} \otimes \wedge^{1} \mathbb{R}^{2 m}\right) \quad \text { for } k \in\{0, \ldots, m-1\}.
	\end{aligned}
\end{equation*}
Moreover,  the first order potential $V_{0,i}$ has the decomposition $V_{0,i}=d \eta_i+F_i$ with
\begin{equation*}
	\eta_i \in W^{2-m, 2}\left(B^{2 m}, s o(n)\right), \quad F_i \in W^{2-m, \frac{2 m}{m+1}, 1}\left(B^{2 m}, \mathbb{R}^{n \times n} \otimes \wedge^{1} \mathbb{R}^{2 m}\right).
\end{equation*}
Assume that all of coefficients $V_{l,i},w_{l,i},\eta_i,F_i$ converge weakly to $V_l,w_l,\eta,F$ in the respective defining spaces. Then there exists a subsequence of $u_i$ which converges weakly in $W^{m,2}$ to a solution $u$ of \eqref{eq:Longue-Gastel system}.
\end{theorem}

Our notations are standard.
By $A\lesssim B$ we mean there exists a universal constant $C>0$ such that $A\le CB$.

\section{Proofs of the main results}

\subsection{Proof of Corollary \ref{coro:fourth case} }

In this section, we shall give a complete proof of Corollary \ref{coro:fourth case}. Since there is only one technical idea  different from that of \cite{Lamm-Riviere-2008},  our proof is almost verbatim  to \cite[Proof of Theorem 1.5]{Lamm-Riviere-2008}. In the key point we will point out the difference between our proof and theirs.  We encourage the readers to  have a closer look at the proof there for a comparison.

\begin{proof}[Proof of Corollary \ref{coro:fourth case} ]
	
Write 
\[
\ta=\|V\|_{W^{1,2}(B^{4})}+\|w\|_{L^{2}(B^{4})}+\|\om\|_{L^{2}(B^{4})}+\|F\|_{L^{4/3,1}(B^{4})}.
\]
To ease our notation, we shall often omit the defining domain in the various norms. For instance, the notation $W^{2,2}$ without a precise defining domain always means on the entire defining domain of the referred objects.

\textbf{Step 1.} Find a suitable Gauge transform 

By \cite[Equation (2.3)]{Lamm-Riviere-2008},
there exists $\Om\in W^{1,2}(B^{4},so_{m}\otimes\wedge^{1}\R^{4})$
such that
\begin{equation}
\begin{cases}
d^{\ast}\Om=\om & \text{in }B^{4},\\
\|\Om\|_{W^{1,2}}\le c\|\om\|_{L^{2}}.
\end{cases}\label{eq: 2.3}
\end{equation}
Then an application of \cite[Theorem A.5]{Lamm-Riviere-2008} gives that, for $\ep_{m}$
sufficiently small with $\|\om\|_{L^{2}}\le\ep_{m}$, there exist
$U\in W^{2,2}(B^{4},so_{m})$ and $P=e^{U}\in W^{2,2}(B^{4},SO_{m})$
and $\xi\in W^{2,2}(B^{4},so_{m}\otimes\wedge^{2}\R^{4})$ and $c_{m}>0$
such that 
\[
\begin{cases}
\Om=P^{-1}dP+P^{-1}d^{\ast}\xi P & \text{in }B^{4},\\
d\left(i_{\pa B^{4}}^{\ast}\ast\xi\right)=0 & \text{on }\pa B^{4},\\
\|P\|_{W^{2,2}}+\|\xi\|_{W^{2,2}}\le c_{m}\|\Om\|_{W^{1,2}}\le C_{m}\ep_{m}.
\end{cases}
\]

\textbf{Step 2.} Rewrite $\om$ and $W$ 

Direct computation gives\footnote{Note that $d^{\ast}\al=-div(\al)$ for any 1-form $\al$. Hence $d^{\ast}(f\al)=fd^{\ast}\al-df\cdot\al$
for any function $f$ and 1-form $\al$. } 
\[
\begin{aligned}\om & =d^{\ast}\Om=d^{\ast}\left(P^{-1}dP+P^{-1}d^{\ast}\xi P\right).\\
 & =-P^{-1}\De P-\langle dP^{-1},dP\rangle-\langle P^{-1}d^{\ast}\xi,dP\rangle-\langle dP^{-1},d^{\ast}\xi P\rangle.
\end{aligned}
\]
Thus 
\begin{equation}
W=d\om+F=-P^{-1}d\De P+K_{1},\label{eq: 2.5}
\end{equation}
where 
\[
K_{1}=-dP^{-1}\De P-d\langle dP^{-1},dP\rangle-d\langle P^{-1}d^{\ast}\xi,dP\rangle-d\langle dP^{-1},d^{\ast}\xi P\rangle+F.
\]
Similar to \cite[Equation (2.9)]{Lamm-Riviere-2008}, using the improved Sobolev embedding 
\[
W^{1,2}(B^{4})\subset L^{4,2}(B^{4}),
\]
one can easily verify that $K_{1}\in L^{4/3,1}(B^{4})$ with 
\[
\|K_{1}\|_{L^{4/3,1}(B^{4})}\le c\|\om\|_{L^{2}(B^{4})}+c\|F\|_{L^{4/3,1}(B^{4})}\le c_{m}\ta.
\]

\textbf{Step 3.} Reduce to an equivalent problem 

Suppose now $A,B$ solves \eqref{eq:fourth order for CL}. Let $\tilde{A}=AP^{-1}$. Then, using (\ref{eq: 2.5}),
$\tilde{A}$ and $B$ solves the equivalent equation
\begin{equation}
d\De\tilde{A}+\De\tilde{A}K_{2}+\na^{2}\tilde{A}K_{3}+d\tilde{A}K_{4}+\tilde{A}K_{5}=d^{\ast}BP^{-1}\label{eq: 2.8}
\end{equation}
in $B^{4}$, where
\begin{equation}
\|K_{2}\|_{W^{1,2}(B^{4})}+\|K_{3}\|_{W^{1,2}(B^{4})}+\|K_{4}\|_{L^{2}(B^{4})}+\|K_{5}\|_{L^{4/3,1}(B^{4})}<c\ta.\label{eq: 2.9}
\end{equation}
Note that by (\ref{eq: 2.8}),
\begin{equation}
d^{\ast}B=\left(d\De\tilde{A}+\De\tilde{A}K_{2}+\na^{2}\tilde{A}K_{3}+d\tilde{A}K_{4}+\tilde{A}K_{5}\right)P.\label{eq: 2.8-2}
\end{equation}

At this moment, different from the cut-off argument of \cite[last paragraph on page 251]{Lamm-Riviere-2008},
we shall use an extension argument as follows. Firstly, we extend $K_{i}$ ($2\le i\le5$)
into $\R^{4}$ with compact support in $B_{2}^{4}$ in a norm-bounded way, for which we denoted by $\tilde{K}_{i}$, such that 
\begin{equation}
\|\tilde{K}_{2}\|_{W^{1,2}(B_{2}^{4})}+\|\tilde{K}_{3}\|_{W^{1,2}(B_{2}^{4})}+\|\tilde{K}_{4}\|_{L^{2}(B_{2}^{4})}+\|\tilde{K}_{5}\|_{L^{4/3,1}(B_{2}^{4})}<c\ta.\label{eq: 2.9-2}
\end{equation}
We next let $\tilde{U}$ be a norm-bounded extension of $U$ into $B_{2}^{4}$ with compact
set so that $\tilde{U}\in W^{2,2}(B_2^4,so_m)$. Then $\tilde{P}=e^{\tilde{U}}\in W^{2,2}(B_{2}^{4},SO_{m})$ with norm bounded by the corresponding norm of $P$. In order to find a solution $(\tilde{A},B)$ for \eqref{eq: 2.8}-\eqref{eq: 2.8-2}, we solve the corresponding problem in the extended region $B_2^4$ as follows:
\begin{equation}
d\De\tilde{A}+\De\tilde{A}\tilde{K}_{2}+\na^{2}\tilde{A}\tilde{K}_{3}+d\tilde{A}\tilde{K}_{4}+\tilde{A}\tilde{K}_{5}=d^{\ast}B\tilde{P}^{-1}\qquad\text{in }B_{2}^{4}.\label{eq: 2.8-1}
\end{equation}
It is clear that each solution of \eqref{eq: 2.8-1} automatically  solves \eqref{eq: 2.8} in $B^{4}$. 

To solve \eqref{eq: 2.8-1}, we turn to look for $(\bar{A},B)$ such
that
\[
\begin{cases}
\tilde{A}=\bar{A}+id\\
dB=0
\end{cases}\qquad\text{in }B_{2}^{4}.
\]
Then $(\bar{A},B)$ necessarily solves the equivalent system 
\begin{equation}
\begin{cases}
\De^{2}\bar{A}=d^{\ast}\left(\De\bar{A}\tilde{K}_{2}+\na^{2}\bar{A}\tilde{K}_{3}+d\bar{A}\tilde{K}_{4}+\bar{A}\tilde{K}_{5}+\tilde{K}_{5}-d^{\ast}B\tilde{P}^{-1}\right) & \text{in }B_{2}^{4},\\
\De B=d\left[\left(d\De\bar{A}+\De\bar{A}\tilde{K}_{2}+\na^{2}\bar{A}\tilde{K}_{3}+d\bar{A}\tilde{K}_{4}+\bar{A}\tilde{K}_{5}+\tilde{K}_{5}\right)\tilde{P}\right] & \text{in }B_{2}^{4},\\
dB=0 & \text{in }B_{2}^{4}.
\end{cases}\label{eq: 2.11-1}
\end{equation}
To solve this system, we impose the following boundary value assumptions
\begin{equation}
\begin{cases}
\bar{A}=\frac{\pa\De\bar{A}}{\pa\nu}=0 & \text{on }\pa B_{2}^{4},\\
\int_{B_{2}^{4}}\De\bar{A}=0,\\
i_{\pa B_{2}^{4}}^{\ast}(\ast B)=0.
\end{cases}\label{eq: 2.11-2}
\end{equation}

\textbf{Step 4.} Solve the equivalent system \eqref{eq: 2.11-1} with boundary value \eqref{eq: 2.11-2}

We will use the fixed point argument as in \cite{Lamm-Riviere-2008} to find a solution $(\bar{A},B)$
for problem \eqref{eq: 2.11-1}-\eqref{eq: 2.11-2}. To this end, we introduce the Banach space $\mathbb{H}=(\mathbb{H},\|\cdot\|_{\mathbb{H}})$ as follows:
\[
\mathbb{H}=\left\{ (u,v)\in W^{2,2}\cap L^{\wq}(B_{2}^{4},M_{m})\times W^{1,4/3}(B_{2}^{4},M_{m}\otimes\wedge^{2}\R^{4}):(u,v)\text{ satisfies }\eqref{eq: 2.11-2}\right\} 
\]
with norm given by
\[
\|(u,v)\|_{\H}\equiv\|u\|_{W^{2,2}(B^{4})}+\|u\|_{L^{\wq}(B^{4})}+\|v\|_{W^{1,4/3}(B^{4})}.
\]
Then, for any $(u,v)\in\H$, there exists a unique solution $\bar{u}\in W^{2,2}(B^{4})$
satisfying 
\[
\begin{cases}
\De^{2}\bar{u}=d^{\ast}\left(\De u\tilde{K}_{2}+\na^{2}u\tilde{K}_{3}+du\tilde{K}_{4}+u\tilde{K}_{5}-d^{\ast}v\tilde{P}^{-1}+\tilde{K}_{5}\right) & \text{in }B_{2}^{4},\\
\bar{u}=\frac{\pa\De\bar{u}}{\pa\nu}=0 & \text{on }\pa B_{2}^{4},\\
\int_{B_{2}^{4}}\De\bar{u}=0.
\end{cases}
\]
Note that $u\in W^{2,2}$ and $v\in W^{1,4/3}$ implies that 
\[
f\equiv\De u\tilde{K}_{2}+\na^{2}u\tilde{K}_{3}+du\tilde{K}_{4}+u\tilde{K}_{5}+\tilde{K}_{5}\in L^{4/3,1}(B_{2}^{4})
\]
with 
\[
\|f\|_{L^{4/3,1}}\le c\ta\left(\|u\|_{W^{2,2}}+\|u\|_{L^{\wq}}+1\right).
\]
and 
\[
d^{\ast}(d^{\ast}v\tilde{P}^{-1})=-\langle d^{\ast}v,d\tilde{P}^{-1}\rangle=\pm\ast(d\ast v\wedge d\tilde{P}^{-1}).
\]
Hence, applying \cite[Lemma A.3]{Lamm-Riviere-2008} (with $w=\tilde{P}^{-1}$,
$p=4/3$, $q=4$), we deduce that $\bar{u}\in W^{3,(4/3,1)}(B_{2}^{4})$
with 
\[
\|\bar{u}\|_{L^{\wq}}+\|\bar{u}\|_{W^{2,2}}+\|\De\bar{u}\|_{L^{2,1}}+\|d\De\bar{u}\|_{L^{\frac{4}{3},1}}\le c\ta(\|u\|_{W^{2,2}}+\|u\|_{L^{\wq}}+1+\|v\|_{W^{1,\frac{4}{3}}}).
\]
Remind that all the norms in the above and below are taken over $B_{2}^{4}$.
Then \cite[Lemma A.1]{Lamm-Riviere-2008} implies that the equations 
\[
\begin{cases}
\De\bar{v}=d\left[\left(d\De\bar{u}+\De\bar{u}\tilde{K}_{2}+\na^{2}\bar{u}\tilde{K}_{3}+d\bar{u}\tilde{K}_{4}+\bar{u}\tilde{K}_{5}+\tilde{K}_{5}\right)\tilde{P}\right] & \text{in }B_{2}^{4},\\
d\bar{v}=0 & \text{in }B_{2}^{4},\\
i_{\pa B_{2}^{4}}^{\ast}(\ast\bar{v})=0
\end{cases}
\]
has a unique solution $\bar{v}\in W^{1,4/3}(B_{2}^{4},M_{m}\otimes\wedge^{2}\R^{4})$
satisfying 
\[
\|d\bar{v}\|_{L^{4/3,1}}\le c\left(\|\bar{u}\|_{W^{3,\frac{4}{3},1}}+\ta\right)\le c\ta\left(\|u\|_{W^{2,2}}+\|u\|_{L^{\wq}}+\|v\|_{W^{1,\frac{4}{3}}}+1\right).
\]
Hence, for any $(u,v)\in\H$, there exists a unique $(\bar{u},\bar{v})\in\H$
solves
\[
\begin{cases}
\De^{2}\bar{u}=d^{\ast}\left(\De u\tilde{K}_{2}+\na^{2}u\tilde{K}_{3}+du\tilde{K}_{4}+u\tilde{K}_{5}-d^{\ast}v\tilde{P}^{-1}+\tilde{K}_{5}\right) & \text{in }B_{2}^{4},\\
\De\bar{v}=d\left[\left(d\De\bar{u}+\De\bar{u}\tilde{K}_{2}+\na^{2}\bar{u}\tilde{K}_{3}+d\bar{u}\tilde{K}_{4}+\bar{u}\tilde{K}_{5}+\tilde{K}_{5}\right)\tilde{P}\right] & \text{in }B_{2}^{4},\\
d\bar{v}=0 & \text{in }B_{2}^{4},
\end{cases}
\]
with boundary value given by \eqref{eq: 2.11-2}. Moreover, there
exists $c>0$, such that 
\begin{equation}
\|(\bar{u},\bar{v})\|_{\H}\le c\ta(\|(u,v)\|_{\H}+1).\label{eq: apriori estimate}
\end{equation}

Now set 
\[
\X=\left\{ (u,v)\in\H:\|(u,v)\|_{\H}\le1\right\} 
\]
and define the mapping $T\colon \H\to\H$ by 
\[
T(u,v)=(\bar{u},\bar{v}).
\]
The apriori estimate \eqref{eq: apriori estimate} implies that we
can choose $\ep_{m}\ll1$ such that $c\ta\le1/2$, and so $T(\X)\subset\X$.
Similar argument implies that $T$ is a contraction operator on $\X$,
and 
\[
\|T(u_{1},v_{1})-T(u_{2},v_{2})\|_{\H}\le\frac{1}{2}\|(u_{1},v_{1})-(u_{2},v_{2})\|_{\H}.
\]
Therefore $T$ is a contraction mapping on $\X$. The standard fixed
point theorem gives a unique $(\bar{A},B)\in\X$ with
\[
T(\bar{A},B)=(\bar{A},B).
\]
This solves \eqref{eq: 2.11-1}-\eqref{eq: 2.11-2}. Moreover,
by \eqref{eq: apriori estimate}, we get $\bar{A}\in W^{3,4/3,1}(B_{2}^{4})$
and $B\in W^{1,4/3,1}(B_{2}^{4})$, and 
\[
\|(\bar{A},B)\|_{\H}\le c\ta\le\ep_{m}.
\]

\textbf{Step 5.} Solve the orginal system

Let $(\bar{A},B)$ be the solution founded in \textbf{Step 4} and set 
\[
\tilde{A}=\bar{A}+id.
\]
Then $\tilde{A}\in W^{3,4/3,1}(B_{2}^{4})$, $B\in W^{1,4/3,1}(B_{2}^{4})$
and they satisfy 
\[
d^{\ast}\left(d\De\tilde{A}+\De\tilde{A}\tilde{K}_{2}+\na^{2}\tilde{A}\tilde{K}_{3}+d\tilde{A}\tilde{K}_{4}+\tilde{A}\tilde{K}_{5}-d^{\ast}B\tilde{P}^{-1}\right)=0\qquad\text{in }B_{2}^{4}.
\]
Thus the nonlinear Hodge decomposition (see \cite[Theorem 2.4.14]{Schwarz-book}) implies that 
\[
d\De\tilde{A}+\De\tilde{A}\tilde{K}_{2}+\na^{2}\tilde{A}\tilde{K}_{3}+d\tilde{A}\tilde{K}_{4}+\tilde{A}\tilde{K}_{5}-d^{\ast}B\tilde{P}^{-1}=d^{\ast}C+h,
\]
where $C\in W^{1,4/3,1}(B_{2}^{4},M_{m}\otimes\wedge^{2}\R^{4})$
satisfies $i_{\pa B_{2}^{4}}^{\ast}(\ast C)=0$ and $h$ is a harmonic
2-form in $B_{2}^{4}$. As $\tilde{K}_{i}\equiv0$ ($2\le i\le5$)
in a $\delta$-neighborhood of $\pa B_{2}^{4}$ , we have
\[
h+d^{\ast}C=d\De\tilde{A}-d^{\ast}B\qquad\text{for } 2-\delta<|x|\le2.
\]
Since $i_{\pa B^{4}}^{\ast}\circ d=d\circ i_{\pa B^{4}}^{\ast}$ and
$i_{\pa B_{2}^{4}}^{\ast}(\ast C)=0$, we deduce 
\[
i_{\pa B_{2}^{4}}^{\ast}\left(\ast h\right)=(d\De\tilde{A}\cdot\nu)\ast vol_{\pa B_{2}^{4}}+di_{\pa B_{2}^{4}}^{\ast}\ast B=0.
\]
This implies that $h\equiv0$, since $B_{2}^{4}$ has trivial homotopy groups, from which we conclude that 
\[
d\De\tilde{A}+\De\tilde{A}\tilde{K}_{2}+\na^{2}\tilde{A}\tilde{K}_{3}+d\tilde{A}\tilde{K}_{4}+\tilde{A}\tilde{K}_{5}-d^{\ast}B\tilde{P}^{-1}=d^{\ast}C
\]
with $i_{\pa B_{2}^{4}}^{\ast}(\ast C)=0$. By the same argument as that of \cite{Lamm-Riviere-2008}, we infer that $C\equiv0$ by choosing $\ep_{m}$ sufficiently small. The proof is complete. 
\end{proof}

\subsection{Proof of Theorem \ref{thm:general even order}}
In this section, we give a sketch of the proof of Theorem \ref{thm:general even order}, since the proof is very similar to \cite[Proof of Theorem 4.1(i)]{deLongueville-Gastel-2019}.

\textbf{Step 1}. Find a suitable Gauge transform.

Given $\eta\in W^{2-m,2}(B^{2m},so(n))$, we extended it in a norm-bounded way to $B^{2m}_2$ with compact support so that its extension $\hat{\eta}\in W^{2-m,2}(B^{2m}_2,so(n))$. Reasoning as in \cite[Second paragraph, page 11]{deLongueville-Gastel-2019}, we can find $\Omega\in W^{m-1,2}(B^{2m}_2,so(n)\wedge^1\R^{2m})$ such that
\begin{equation}
	\begin{cases}
		\Delta^{m-2}d^\ast\Om=-\hat{\eta} & \text{in }B^{2m}_2,\\
		\|\Om\|_{W^{m-1,2}(B^{2m}_2)}\le c\|\hat{\eta}\|_{W^{2-m,2}(B_2^{2m})}.
	\end{cases}\label{eq:1}
\end{equation}
When the smallness assumption \eqref{eq:smallness assumption} is satisfied for some sufficiently small $\ep_{m}$, we may apply \cite[Theorem 2.4]{deLongueville-Gastel-2019} to find $P\in W^{m,2}(B^{2m},SO(n))$ and $\xi\in W^{m,2}(B^{2m},so(n)\otimes \wedge^1\R^{2m})$ such that
$$\Omega=P^{-1}dP+P^{-1}d^\ast \xi P\qquad \text{on}\  B^{2m}.$$
Moreover, $\|dP\|_{W^{m-1,2}(B^{2m})}+\|d^\ast \xi\|_{W^{m-1,2}(B^{2m})}\leq c\|\Omega\|_{W^{m-1,2}(B^{2m}_2)}$.

 \textbf{Step 2}. Rewrite $V_0$.
 
 By definition, we have
 $$V_0=d\eta+F=-d\Delta^{m-2}d^\ast\big(P^{-1}dP+P^{-1}d^\ast \xi P\big)+F\quad \text{on }B^{2m}.$$
 Similar to \cite[Third paragraph, page 11]{deLongueville-Gastel-2019}, we could further write it as
 $$V_0=-P^{-1}d\Delta^{m-1}P+K,$$
 with $K\in W^{2-m,\frac{2m}{m+1},1}(B^{2m})$.
 
 \textbf{Step 3.} Reduce to an equivalent problem.
 
 Suppose $A, B$ solves \eqref{eq:for CL general order}. Let $\tilde{A}=AP^{-1}$. Then $(\tilde{A},B)$ solves the equivalent equation (see \cite[Equation (11)]{deLongueville-Gastel-2019})
 \begin{equation}\label{eq:equivalent system general order}
 	d\Delta^{m-1}\tilde{A}+\sum_{j=0}^{2m-2}\langle D^j\tilde{A},K_j\rangle+K_0=d^\ast BP^{-1},
 \end{equation}
 with all the coefficent functions $K_0,\cdots,K_{m-2}$ bounded by $c\theta$ in the respective norm spaces 
 $$K_0\in W^{2-m,\frac{2m}{m+1},1}(B^{2m}),K_j\in W^{j+1-m,2}(B^{2m})\quad \text{for }j\in \{1,\cdots,m-2\}.$$
 
 Different from the cut-off argument used in \cite[Last paragraph, page 11]{deLongueville-Gastel-2019}, we again use an extension argument similar to the proof of Corollary \ref{coro:fourth case}. Namely, we extend all the coefficient functions $K_j$, $j=0,\cdots,m-2$ into $\R^{2m}$ with compact support in $B_2^{2m}$ in a norm-bounded way, for which we denoted by $\tilde{K}_j$. Then we use a similar extension $\tilde{P}$ of $P$ such that $\tilde{P}\in W^{m,2}(B_2^{2m},SO(n))$. 
 
 In order to find a solution $(\tilde{A},B)$ for \eqref{eq:equivalent system general order}, we solve the corresponding problem in the extended region $B_2^{2m}$ as follows:
 \begin{equation}
 	d\Delta^{m-1}\tilde{A}+\sum_{j=0}^{2m-2}\langle D^j\tilde{A},\tilde{K}_j\rangle+\tilde{K}_0=d^\ast B\tilde{P}^{-1}\qquad\text{in }B_{2}^{2m}.\label{eq:2}
 \end{equation}
 It is clear that each solution of \eqref{eq:2} automatically  solves \eqref{eq:equivalent system general order} in $B^{2m}$.
 
 From now on, the remaining steps are identical to that used in \cite{deLongueville-Gastel-2019}. One uses a fix point argument to find solutions $(\tilde{A},B)$ that solves \eqref{eq:2}; for details, see \cite[Proof of Lemma 4.2]{deLongueville-Gastel-2019}. The proof is thus complete.

\subsection{Proof of Theorem \ref{thm:weak compactness} }
In this section, we shall give a proof of Theorem \ref{thm:weak compactness}, following closely the approach of Rivi\`ere \cite{Riviere-2007} for the second order system \eqref{eq:Riviere 2007}. 

\textbf{Step 1}. Cover the ball $B^{2m}$ with small ``good balls'' $B_j^{2m}\subset B^{2m}$.

We may cover the ball $B^{2m}$ by balls $B_j^{2m}\subset B^{2m}$ in such a way that in every $B_j^{2m}$ we may assume that
$$\theta_{B_j^{2m}}<\ep_{m}.$$
This collection of balls will cover the whole $B^{2m}$ outside a finite set of points.

\textbf{Step 2}. Derive the convergence on each good ball $B_j:=B_j^{2m}$.

Since $\theta_{B_j^{2m}}<\ep_{m}$, we may apply Theorem \ref{thm:general even order} to find  $A_i=A_{i,j}\in W^{m,2}\cap L^\infty(B_{j}^{2m},Gl(n))$ and $B_i=B_{i,j}\in W^{2-m,2}(B^{2m}_{j},\R^{n\times n}\otimes \wedge^2\R^{2m})$  such that
\begin{equation}\label{eq:for CL general order seq}
	\Delta^{m-1}dA_i+\sum_{k=0}^{m-1}(\Delta^k A_i)V_{k,i}-\sum_{k=0}^{m-2}(\Delta^k dA_i)w_{k,i}=\delta B_i\qquad \text{in }B_j.
\end{equation}
Moreover, $u_i$ is a distributional solution of 
\begin{equation}\label{eq:conservation law of D-G seq}
	\begin{aligned}
		-A_if_i&=\delta\Big[\sum_{l=0}^{m-1}\left(\Delta^{l} A_i\right) \Delta^{m-l-1} d u_i-\sum_{l=0}^{m-2}\left(d \Delta^{l} A_i\right) \Delta^{m-l-1} u_i \\ &\qquad -\sum_{k=0}^{m-1} \sum_{l=0}^{k-1}\left(\Delta^{l} A_i\right) \Delta^{k-l-1} d\left\langle V_{k,i}, d u_i\right\rangle+\sum_{k=0}^{m-1} \sum_{l=0}^{k-1}\left(d \Delta^{l} A_i\right) \Delta^{k-l-1}\left\langle V_{k,i}, d u_i\right\rangle \\ &\qquad -\sum_{k=0}^{m-2} \sum_{l=0}^{k-2}\left(\Delta^{l} A_i\right) d \Delta^{k-l-1} \delta\left(w_{k,i} d u_i\right)+\sum_{k=0}^{m-2} \sum_{l=0}^{k-2}\left(d \Delta^{l} A_i\right) \Delta^{k-l-1} \delta\left(w_{k,i} d u_i\right) \\ &\qquad -\langle B_i, d u_i\rangle\Big].
	\end{aligned}
\end{equation}
Up to a subsquence if necessary, we may assume that $A_i$ converges weakly to $A\in W^{m,2}$ and $B_i$ converges weakly in $W^{3-m,\frac{2m}{m+1}}$.

It is easy to check that 
\begin{equation}\label{eq:LR III 3 and 4}
	\begin{aligned}
		\Delta^{m-1}d A_{i} &\to \Delta^{m-1}d A \quad \text{in } \mathcal{D}'(B_j),\\
		(\Delta^k A_i)V_{k,i}&\to (\Delta^k A)V_{k}\quad \text{in } \mathcal{D}'(B_j),\\
		(\Delta^k dA_i)w_{k,i}&\to (\Delta^k dA)w_{k}\quad \text{in } \mathcal{D}'(B_j),\\
		\delta(B_{k,i}) &\to \delta(B_i) \quad \text{in } \mathcal{D}'(B_j).
	\end{aligned}
\end{equation}
 Similarly, we have
\begin{equation}\label{eq:LR III 4}
	\begin{aligned}
		\sum_{l=0}^{m-1}\left(\Delta^{l} A_i\right) \Delta^{m-l-1} d u_i&\to \sum_{l=0}^{m-1}\left(\Delta^{l} A\right) \Delta^{m-l-1} d u \quad \text{in } \mathcal{D}'(B_j),\\
		\sum_{l=0}^{m-2}\left(d \Delta^{l} A_i\right) \Delta^{m-l-1} u_i&\to \sum_{l=0}^{m-2}\left(d \Delta^{l} A\right) \Delta^{m-l-1} u\quad \text{in } \mathcal{D}'(B_j),\\
		\sum_{k=0}^{m-1} \sum_{l=0}^{k-1}\left(\Delta^{l} A_i\right) \Delta^{k-l-1} d\left\langle V_{k,i}, d u_i\right\rangle&\to\sum_{k=0}^{m-1} \sum_{l=0}^{k-1}\left(\Delta^{l} A\right) \Delta^{k-l-1} d\left\langle V_{k}, d u\right\rangle\quad \text{in } \mathcal{D}'(B_j),\\
		\sum_{k=0}^{m-1} \sum_{l=0}^{k-1}\left(d \Delta^{l} A_i\right) \Delta^{k-l-1}\left\langle V_{k,i}, d u_i\right\rangle&\to \sum_{k=0}^{m-1} \sum_{l=0}^{k-1}\left(d \Delta^{l} A\right) \Delta^{k-l-1}\left\langle V_{k}, d u\right\rangle\quad \text{in } \mathcal{D}'(B_j)\\
		\sum_{k=0}^{m-2} \sum_{l=0}^{k-2}\left(\Delta^{l} A_i\right) d \Delta^{k-l-1} \delta\left(w_{k,i} d u_i\right)&\to\sum_{k=0}^{m-2} \sum_{l=0}^{k-2}\left(\Delta^{l} A\right) d \Delta^{k-l-1} \delta\left(w_{k} d u\right)\quad \text{in } \mathcal{D}'(B_j),\\
		\sum_{k=0}^{m-2} \sum_{l=0}^{k-2}\left(d \Delta^{l} A_i\right) \Delta^{k-l-1} \delta\left(w_{k,i} d u_i\right)&\to \sum_{k=0}^{m-2} \sum_{l=0}^{k-2}\left(d \Delta^{l} A\right) \Delta^{k-l-1} \delta\left(w_{k} d u\right)\quad \text{in } \mathcal{D}'(B_j)\\
		\langle B_i, d u_i\rangle &\to \langle B, d u\rangle \quad \text{in } \mathcal{D}'(B_j)
	\end{aligned}
\end{equation}
and
\begin{equation}\label{eq:LR III 5}
	A_{i}f_i\to 0 \quad \text{in } \mathcal{D}'(B_j).
\end{equation}
Combining \eqref{eq:for CL general order seq} - \eqref{eq:LR III 5} all together, we infer from Theorem \ref{thm:general even order} that $u$ is a solution to \eqref{eq:Longue-Gastel system} in $B_j$.

\textbf{Step 3}. Remove the singular set on $B^{2m}$.

By \textbf{Step 2}, we know the distribution 
$$T:=\Delta^{m}u-\sum_{l=0}^{m-1}\Delta^{l}\left\langle V_{l},du\right\rangle -\sum_{l=0}^{m-2}\Delta^{l}\delta\left(w_{l}du\right)$$
vanishes on each good ball $B_j$ and hence $T$ only supports at finitely many points in $B^{2m}$. On the other hand, we know that $T\in H^{-m}+L^1$ and so it has to be identically zero on $B^{2m}$. This completes the proof.

\textbf{Acknowledgments}. The authors are grateful to Prof.~T. Rivi\`ere and Prof.~T. Lamm for bringing \textbf{Problem A} to their atttension and for several helpful communications during the preparation of this work.

\end{document}